\documentclass[11pt, reqno]{amsart}

\usepackage[
  letterpaper,
  left=1.2in,
  right=1.2in,
  top=1.1in,
  bottom=1.1in
]{geometry}

\usepackage{cite}
\usepackage{etoolbox}
\usepackage{amsmath,amssymb,graphicx,float}
\usepackage{algorithm}
\usepackage{algpseudocode}
\usepackage{hyperref}
\floatstyle{ruled}
\newfloat{algorithm}{tbp}{loa}
\floatname{algorithm}{Algorithm}

\hypersetup{
	colorlinks   = true, 
	urlcolor     = black, 
	linkcolor    = black, 
	citecolor   = black 
}

\makeatletter
\patchcmd{\thebibliography}
  {\usecounter{enumiv}}
  {\usecounter{enumiv}\setlength{\itemsep}{.6em}}
  {}{}
\makeatother

\newtheorem{theorem}{Theorem}[section]
\newtheorem{lemma}[theorem]{Lemma}
\newtheorem{proposition}[theorem]{Proposition}
\newtheorem{corollary}[theorem]{Corollary}

\theoremstyle{definition}

\theoremstyle{remark}
\newtheorem{remark}{Remark}

\theoremstyle{plain}

\numberwithin{equation}{section}

\newcommand{\Lstab}{L_{\mathrm{stab}}}

\begin{document}
\raggedbottom

\title[Nonlocal Wound Healing on Bounded Domains]{Stability analysis and critical thresholds in a nonlocal advection diffusion model of wound healing on bounded domains}

\author{Ilyas Mustapha and Mikayla Gregory}
\address{Department of Mathematics, Washington University in St.\,Louis, St.\,Louis, MO, USA}






\begin{abstract}
Nonlocal advection diffusion models describe how wound healing depends on the strength of cell-cell adhesion and the cell sensing radius. Depending on these parameters, the wound may close completely or remain in a persistent state with a spatially nonuniform cell density. The model combines random cell movement, cell proliferation, and cell-cell adhesion within a finite sensing radius. We first derive the dispersion relation on the whole line and use it to compute the critical sensing radius and the corresponding wavelength. This relation shows that diffusion and cell proliferation stabilize the healed state, while cell-cell adhesion can make it unstable. On a bounded interval, the cell density is fixed at the healthy value on an exterior collar. In this setting, we establish conditions that keep the cell density within a positive bounded range, derive a sufficient condition for nonlinear exponential stability that depends on domain size, and formulate and numerically solve the associated linear spectral problem. The analysis shows how reducing the domain size can stabilize the healed state at parameter values for which the whole line model is linearly unstable. Finally, we compute critical sensing radii and domain sizes numerically and compare the sufficient stability bounds with the transitions observed in wound simulations involving perturbations of finite amplitude.

\par\vspace{0.25cm}
\end{abstract}

\maketitle

\fontsize{10}{13}\selectfont

\section{Introduction}
In this work, we use a nonlocal advection diffusion model on bounded domains to study how cell adhesion, cell sensing radius, and domain size affect the stability of the healed state and the conditions for complete healing. The model combines random cell movement, cell proliferation, and cell adhesion to describe wound closure. Experimentally, this process can be studied through scratch assays, in which part of a cell monolayer is removed and the closure of the gap is recorded over time \cite{jonkman2014introduction}. Earlier mathematical models described wound closure by combining diffusion for cell movement with cell proliferation and chemical regulation \cite{sherratt1991mathematical,sherratt1990models}. Other models included tissue remodeling, wound contraction, and angiogenesis to account for further aspects of the healing process \cite{sherratt2002theoretical}. The movement of cells into the wound has also been studied using Fisher-type models, which describe cell fronts that approach a constant speed after an initial transient \cite{maini2004traveling}. Mechanical models of collective cell migration have shown that closure time depends on wound shape and not on area alone \cite{arciero2011continuum,arciero2013using}. Cell adhesion can support this collective movement, but its effect on wound closure depends on how it interacts with random cell movement and proliferation, as well as the parameter values \cite{khain2007role,nardini2016modeling}.

Armstrong, Painter, and Sherratt \cite{armstrong2006continuum} developed a nonlocal continuum model of cell adhesion in which the velocity, due to adhesion, depends on the surrounding cell density within a finite sensing radius. This nonlocal description was later connected to individual cell movement through spatial jump processes \cite{buttenschon2018space}. The same framework has been used to study cancer invasion \cite{gerisch2008mathematical}, somite formation \cite{armstrong2009adding}, attraction and repulsion between cells \cite{painter2015nonlocal}, and cell sorting and crowding \cite{murakawa2015continuous,carrillo2019population}. The mathematical analysis has also developed alongside these applications. Results on boundedness, existence, and long time behavior have been established for nonlocal adhesion models on bounded and unbounded domains \cite{sherratt2009boundedness,dyson2010existence,dyson2013non}. Further work has provided global solvability and explicit bounds \cite{hillen2018global}, while steady states and bifurcations have been studied to describe changes in the resulting cell distributions \cite{buttenschon2021non}. These developments have made nonlocal adhesion models useful for connecting short-range cell interactions with larger scale patterns of biological movement \cite{painter2024biological}.

Recent models have included interactions among fibroblasts, macrophages, extracellular matrix, and growth factors to study normal and abnormal wound healing, together with questions of stability and local existence of solutions \cite{adebayo2023mathematical,adebayo2024analytical}. Cell adhesion has also been incorporated into models describing the role of keratinocyte growth factor in wound healing \cite{endes2025cellular}. Comparisons with scratch assay data have shown that models with and without adhesion can fit the same data similarly well while giving different estimates of the diffusion coefficient \cite{atici2025effect}. This shows that the interpretation of experimental data can depend on which processes are included in the model. Other work has connected nonlocal cell adhesion to receptor binding and studied the existence and uniqueness of solutions \cite{zhigun2024modelling,rajendran2025local}. Related mathematical studies have examined local approximations and the formation and stability of cell clusters \cite{falco2025nonlocal,buttenschon2024cells}, positivity and global existence \cite{giunta2023positivity}, and pattern formation in models with nonlocal advection \cite{shi2026bifurcation}. Here, we use a scalar model with prescribed parameter values to study how cell adhesion, the sensing radius, and domain size affect healing.

Following the adhesion formulation in \cite{armstrong2006continuum}, we consider the nonlocal advection diffusion model of wound healing given by
\begin{equation}\label{eq:wound_healing_model}
    n_t = \alpha n_{xx}
    - \partial_x \left[
        n(x,t)\int_{-\rho}^{\rho}
        g(n(x+\xi,t))h(\xi)\,d\xi
    \right]
    + f(n).
\end{equation}
Here, $n(x,t)$ is the cell density, $\alpha>0$ is the diffusion coefficient, $\rho>0$ is the cell sensing radius, and $h$ is an odd interaction kernel. The three terms on the right represent random cell movement, cell adhesion, and cell proliferation, respectively. The adhesion integral accounts for interactions with surrounding cells within the sensing radius. In our simulations, we take 
$$f(n)=n(1-n), \text{ and } g(n)=n(\lambda-n),$$ where $\lambda$ controls adhesion and the uniform density $n\equiv1$ represents the healed state. The equation with cell proliferation was analyzed in \cite{dyson2010existence}, where wound closure simulations were also presented. Later work examined how the sensing radius and adhesion strength affect wound closure \cite{webb2022force}, followed by a bifurcation analysis distinguishing uniform from nonuniform
closure \cite{webb2023bifurcation}, in which a critical value of the adhesion strength separates
stable from unstable uniform closure. The present work fixes the adhesion
strength and treats the sensing radius and the domain half-length as the
parameters. The domain half-length has no whole line counterpart and requires the bounded
linearized operator, whose eigenfunctions are not in general the Dirichlet
sine modes.


Nonlocal adhesion models on bounded domains have also been studied in \cite{hillen2020nonlocal,buttenschon2021non}. In \cite{hillen2020nonlocal}, adhesive, repulsive, and no-flux boundary conditions are derived for the nonlocal operator and the resulting problems
are shown to be well posed. We instead prescribe the exterior values at the healed density, which
represents the intact monolayer surrounding a scratch, and use this formulation to compute the spectrum of the bounded linearized operator and the thresholds it determines.

In our setting, we take $\Omega=(-L,L)$ and fix the cell density at $n=1$ on an exterior collar of width $\rho$. This collar provides the density values needed when the adhesion integral samples points outside the interval and represents surrounding cells maintained at the healthy density. Here, $L$ is half the length of the domain and should be distinguished from half the width of the initial wound. We keep the initial wound width fixed and vary $L$ to study how the distance between the wound and the healthy exterior affects healing. The bounded domain results therefore apply to this prescribed exterior condition and need not hold for other boundary conditions.

The whole line analysis provides a reference for the bounded domain problem. The dispersion relation describes whether small perturbations of the healed state grow or decay and gives the first critical sensing radius $\rho^\infty_c(\lambda)$. On a bounded interval, the exterior condition changes the linear stability problem, and the Dirichlet sine functions do not generally diagonalize the nonlocal operator. Simply restricting the whole line dispersion relation to discrete wavenumbers therefore does not give the exact bounded domain threshold. We use the full bounded domain linearized operator to determine the critical sensing radius $\rho^L_c(\lambda)$ and, where it exists, the critical domain half length $L_c(\rho,\lambda)$ at which the healed state loses strict linear stability. An energy estimate using the Poincar\'e inequality gives a sufficient bound $\Lstab(\rho, \lambda)$ below which solutions converge exponentially to the healed state, provided the assumed bounds on cell density hold. This bound depends on the model parameters and the density range, so it does not give a single domain size that guarantees healing for every parameter choice. We compare these predictions with simulations starting from a wound of finite amplitude. The simulations test whether the linear thresholds predict healing for the chosen wound profiles and whether complete healing occurs for domain sizes beyond the sufficient bound $\Lstab$.

The observation time is important when interpreting these simulations. Patterns that appear persistent may still be changing slowly, as shown for a related aggregation model in \cite{potts2024distinguishing}. A wound that has not healed by the final simulation time cannot therefore be assumed to remain unhealed indefinitely. We compute numerical solutions using Fourier continuation for spatial derivatives, efficient evaluation of the adhesion integral \cite{gerisch2010approximation}, and fourth order Runge-Kutta time integration. The development and analysis of Fourier continuation methods are described in \cite{bruno2010high, lyon2010high}, with applications to nonlocal equations on bounded domains in \cite{mustapha2025fourier}. We refer to these studies for the analysis of the methods and focus here on critical parameter values, healing times, and cell density profiles. 

The paper is organized as follows. In Section~\ref{sec:whole line}, we derive the dispersion relation for the healed state on the whole line and use it to define the critical sensing radius and the associated wavelength.
Section~\ref{sec:bounded-formulation} presents the bounded-domain formulation with a prescribed healthy exterior, an invariant density interval, a domain-dependent sufficient condition for nonlinear exponential stability, and the bounded linearized operator that determines the critical sensing radius and length
thresholds. In Section~\ref{sec:numerical-experiments}, we discuss the numerical results and examine how the sensing radius and the domain size affect healing. In Section~\ref{sec:conclusion}, we provide a summary of our approach and outline directions for future research.

\section{Stability analysis on the whole line}
\label{sec:whole line}	
We study the stability of the healed state $n\equiv1$ for the model in \eqref{eq:wound_healing_model}, with $f(n)=n(1-n)$ and $g(n)=n(\lambda-n)$. We derive the dispersion relation on the whole line and use it to define the critical sensing radius and the associated wavelength.

\subsection{Linear stability on the whole line and the dispersion relation}
Consider model~\eqref{eq:wound_healing_model} on the whole line, where $\alpha>0$, $\rho>0$, and $h\in L^1(-\rho,\rho)$ is a real, odd interaction kernel. Let $\bar n$ be a constant satisfying $f(\bar n)=0$. Since $h$ is odd,
the adhesion velocity vanishes at $n=\bar n$. Thus, $\bar n$ is a constant equilibrium of \eqref{eq:wound_healing_model}.

To obtain the linearized equation, write
\[
n(x,t)=\bar n+\varepsilon w(x,t),
\quad |\varepsilon|\ll1,
\]
where $w$ is a smooth perturbation. Define
\[
V[n](x,t)
=
\int_{-\rho}^{\rho}
g\bigl(n(x+\xi,t)\bigr)h(\xi)\,d\xi,
\quad
W(x,t)
=
\int_{-\rho}^{\rho}w(x+\xi,t)h(\xi)\,d\xi.
\]
At the constant equilibrium,
\[
V[\bar n]=0, \quad
\left.\frac{d}{d\varepsilon}
V[\bar n+\varepsilon w]\right|_{\varepsilon=0}
=
g'(\bar n)W.
\]
The derivative of the adhesion flux is consequently
\[
\left.\frac{d}{d\varepsilon}
\bigl[(\bar n+\varepsilon w)
V[\bar n+\varepsilon w]\bigr]
\right|_{\varepsilon=0}
=
wV[\bar n]+\bar n g'(\bar n)W
=
\bar n g'(\bar n)W.
\]
Hence, differentiating \eqref{eq:wound_healing_model} with respect to $\varepsilon$ at $\varepsilon=0$ gives
\begin{equation}
\label{eq:wl-linearized}
w_t
=
\alpha w_{xx}
-
\bar n g'(\bar n)W_x
+
f'(\bar n)w.
\end{equation}

\begin{proposition}[Dispersion relation]
\label{prop:wl-dispersion}
For each $k\in\mathbb R$, the linearized equation
\eqref{eq:wl-linearized} admits a mode
$w(x,t)=e^{ikx+\sigma(k)t}$ with growth rate
\begin{equation}
\label{eq:wl-dispersion-general}
\sigma(k)
=
-\alpha k^2+f'(\bar n)
+\bar n g'(\bar n)kS_\rho(k),
\end{equation}
where
\begin{equation}
\label{eq:wl-sine-transform}
S_\rho(k)
=
2\int_0^\rho h(\xi)\sin(k\xi)\,d\xi.
\end{equation}
The growth rate $\sigma(k)$ is real and even in $k$.
\end{proposition}

\begin{proof}
For $w(x,t)=e^{ikx+\sigma t}$,
\[
w_t=\sigma w,
\quad
w_{xx}=-k^2w,
\quad
w(x+\xi,t)=e^{ik\xi}w(x,t).
\]
It follows that
\[
W(x,t)
=
w(x,t)\int_{-\rho}^{\rho}h(\xi)e^{ik\xi}\,d\xi.
\]
Since
\begin{align*}
\int_{-\rho}^{\rho}h(\xi)e^{ik\xi}\,d\xi
&=
\int_{-\rho}^{\rho}h(\xi)\cos(k\xi)\,d\xi
+
i\int_{-\rho}^{\rho}h(\xi)\sin(k\xi)\,d\xi\\
&=
2i\int_0^\rho h(\xi)\sin(k\xi)\,d\xi\\
&=
iS_\rho(k),
\end{align*}
therefore,
\begin{equation*}
W=iS_\rho(k)w \quad
\text{ and } \quad
W_x=iS_\rho(k)(ik)w=-kS_\rho(k)w.
\end{equation*}
Substitution into \eqref{eq:wl-linearized} yields
\[
\sigma w
=
\left[
-\alpha k^2
+\bar n g'(\bar n)kS_\rho(k)
+f'(\bar n)
\right]w,
\]
which proves \eqref{eq:wl-dispersion-general}.
Since $h$ is real, $S_\rho(k)$ is real. Moreover,
$S_\rho(-k)=-S_\rho(k)$, so $kS_\rho(k)$ and
$\sigma(k)$ are even.
\end{proof}

\begin{remark}
For the healed state $\bar n=1$, the choices $f(n)=n(1-n)$ and $g(n)=n(\lambda-n)$ give
$f'(1)=-1$ and $g'(1)=\lambda-2$. Thus \eqref{eq:wl-linearized} becomes
\begin{equation*}
w_t
=
\alpha w_{xx}-w
-(\lambda-2)\partial_x
\left[
\int_{-\rho}^{\rho}w(x+\xi,t)h(\xi)\,d\xi
\right],
\end{equation*}
and the dispersion relation is
\begin{equation}
\label{eq:wl-healed-dispersion}
\sigma(k;\rho,\lambda)
=
-\alpha k^2-1+(\lambda-2)kS_\rho(k).
\end{equation}
\end{remark}

The term $-\alpha k^2$ describes the damping due to random cell
movement, while $-1$ comes from cell proliferation near the healed
state. The adhesion term can increase or decrease the growth rate,
depending on the sign of $(\lambda-2)kS_\rho(k)$. A mode grows only
when this term exceeds $\alpha k^2+1$. Thus any linear instability
of the healed state in this model must come from the adhesion term.
In particular, when $\lambda=2$, the linear adhesion term vanishes
and $\sigma(k;\rho,2)=-\alpha k^2-1<0$ for every $k$.

For fixed $\alpha$, define the largest growth rate by
\begin{equation}
\label{eq:wl-largest-growth}
F(\rho;\lambda)
=
\max_{k\in\mathbb R}\sigma(k;\rho,\lambda)
=
\max_{k\ge0}\sigma(k;\rho,\lambda).
\end{equation}
This maximum exists since $S_\rho$ is continuous and
\begin{equation*}
|S_\rho(k)|
\le
2\int_0^\rho |h(\xi)|\,d\xi
=
\int_{-\rho}^{\rho}|h(\xi)|\,d\xi
=:H_\rho.
\end{equation*}
Hence,
\[
\sigma(k;\rho,\lambda)
\le
-\alpha k^2-1+|\lambda-2||k|H_\rho
\longrightarrow-\infty
\quad\text{as }|k|\longrightarrow\infty.
\]

The sign of $F$ determines linear stability. If $F<0$, every Fourier mode decays at a rate bounded away from zero. If $F>0$, there is an interval of wavenumbers with positive growth rates. At $F=0$, no Fourier mode grows in the linearized equation, but there is no single positive exponential decay rate that applies to all perturbations. The dispersion relation alone does not determine whether solutions of the nonlinear model converge to the healed state.


These statements also apply to square-integrable perturbations by taking the Fourier transform of \eqref{eq:wl-linearized} in $x$, which gives
\[
\partial_t\widehat w(k,t)
=
\sigma(k)\widehat w(k,t).
\]
Solving this equation yields
\begin{equation}
\label{eq:wl-fourier-evolution}
\widehat w(k,t)
=
e^{\sigma(k)t}\widehat w_0(k),
\end{equation}
where $w_0=w(\cdot,0)$. At the healed state, Plancherel's identity therefore gives
\begin{equation*}
\|w(\cdot,t)\|_{L^2(\mathbb R)}
\le
e^{F(\rho;\lambda)t}
\|w_0\|_{L^2(\mathbb R)}.
\end{equation*}
If $F>0$, continuity gives an interval on which $\sigma(k)\ge\gamma>0$. Choosing nonzero initial Fourier data
supported in that interval gives an $L^2$ perturbation growing
at least as fast as $e^{\gamma t}$. Real perturbations can be
obtained by choosing conjugate-symmetric data on this interval
and its reflection.

An upper bound on $S_\rho$ gives a sufficient stability condition
without finding the maximum in \eqref{eq:wl-largest-growth}.

\begin{proposition}[Sufficient linear stability on the whole line]
\label{prop:wl-sufficient-stability}
Let $\bar n$ be a constant equilibrium with $f'(\bar n)<0$, and set
$
H_\rho=\int_{-\rho}^{\rho}|h(\xi)|\,d\xi.
$
If
\begin{equation}
\label{eq:wl-sufficient-condition}
\bigl(\bar n g'(\bar n)\bigr)^2H_\rho^2
<
-4\alpha f'(\bar n),
\end{equation}
then
\[
\sigma(k)\le-\eta<0
\quad\text{for every }k\in\mathbb R,
\]
where
\begin{equation*}
\eta
=
-f'(\bar n)
-
\frac{\bigl(\bar n g'(\bar n)\bigr)^2H_\rho^2}{4\alpha}.
\end{equation*}
Consequently, the solution of the linearized equation
\eqref{eq:wl-linearized} with initial perturbation
$w_0\in L^2(\mathbb R)$ satisfies
\begin{equation*}
\|w(\cdot,t)\|_{L^2(\mathbb R)}
\le
e^{-\eta t}\|w_0\|_{L^2(\mathbb R)},
\quad t\ge0.
\end{equation*}
\end{proposition}

\begin{proof}
Set $a=\bar n g'(\bar n)$ and $B=|a|H_\rho$. Then, for every $k$, Proposition~\ref{prop:wl-dispersion} gives
\begin{align*}
\sigma(k)
&=-\alpha k^2+f'(\bar n)+akS_\rho(k)\\
&\le -\alpha k^2+f'(\bar n)+B|k|\\
&= -\alpha\left(|k|-\frac{B}{2\alpha}\right)^2 + f'(\bar n)+\frac{B^2}{4\alpha}\\
&\le f'(\bar n)+\frac{B^2}{4\alpha}\\
&=-\eta.
\end{align*}
Condition \eqref{eq:wl-sufficient-condition} implies $\eta>0$. Hence, \eqref{eq:wl-fourier-evolution} and Plancherel's identity
give
\[
\begin{aligned}
\|w(\cdot,t)\|_{L^2(\mathbb R)}^2
&=
\int_{\mathbb R}
e^{2\sigma(k)t}|\widehat w_0(k)|^2\,dk
\le
e^{-2\eta t}
\int_{\mathbb R}|\widehat w_0(k)|^2\,dk
=
e^{-2\eta t}\|w_0\|_{L^2(\mathbb R)}^2.
\end{aligned}
\]
\end{proof}

\begin{remark}
    At the healed state, condition \eqref{eq:wl-sufficient-condition} reduces to
    \begin{equation*}
    (\lambda-2)^2H_\rho^2<4\alpha.
    \end{equation*}
\end{remark}

This condition is sufficient, but it is not necessary. Its failure does not establish instability.

\subsection{Critical sensing radius and wavelength.}
For fixed $\alpha$, $\lambda$, and a specified kernel family, define the critical sensing radius over the whole line by
\begin{equation}
\label{eq:wl-critical-radius}
\rho_c^\infty(\lambda) =
\inf\{\rho>0:F(\rho;\lambda)>0\}.
\end{equation}
If this set is empty, we take $\rho_c^\infty(\lambda)=+\infty$.
Any dependence of the kernel on $\rho$, including its
normalization, is included when evaluating $S_\rho$ and $F$.

\begin{lemma}[Growth rate at the critical sensing radius]
\label{lem:wl-critical-growth}
Let $\rho_c^\infty(\lambda)$ be as defined in
\eqref{eq:wl-critical-radius}. Suppose that $0<\rho_c^\infty(\lambda)<\infty$ and that $F(\cdot;\lambda)$ is continuous at $\rho_c^\infty(\lambda)$. Then
\[
F(\rho_c^\infty(\lambda);\lambda)=0.
\]
\end{lemma}

\begin{proof}
Suppose first that
$F(\rho_c^\infty(\lambda);\lambda)>0$.
By continuity, $F(\rho;\lambda)>0$ for some
$0<\rho<\rho_c^\infty(\lambda)$. This contradicts the
definition of $\rho_c^\infty(\lambda)$ as the infimum
of the radii for which $F$ is positive.

If instead $F(\rho_c^\infty(\lambda);\lambda)<0$,
continuity implies that $F(\rho;\lambda)<0$ for all
radii sufficiently close to $\rho_c^\infty(\lambda)$.
However, the definition of the infimum guarantees that
radii with $F(\rho;\lambda)>0$ occur arbitrarily close
to $\rho_c^\infty(\lambda)$ from above. This is again
a contradiction. Therefore, $F(\rho_c^\infty(\lambda);\lambda)=0.$
\end{proof}

Let $k_c\ge0$ be a global maximizer at this critical radius.
Since $\sigma(0;\rho,\lambda)=-1$, it follows that $k_c>0$.
The maximum is therefore attained at an interior wavenumber,
and
\[
\sigma\bigl(k_c;\rho_c^\infty(\lambda),\lambda\bigr)=0,
\quad
\partial_k\sigma\bigl(k_c;\rho_c^\infty(\lambda),\lambda\bigr)=0.
\]
Differentiation under the integral is valid because
$|\xi h(\xi)|$ is integrable on the finite sensing interval.
It gives
\[
\partial_k S_\rho(k)
=
2\int_0^\rho \xi h(\xi)\cos(k\xi)\,d\xi.
\]
Thus the critical radius and wavenumber satisfy
\begin{equation}
\label{eq:wl-critical-pair}
\begin{aligned}
(\lambda-2)k_cS_{\rho_c^\infty}(k_c)
&=\alpha k_c^2+1,\\
(\lambda-2)\left[
S_{\rho_c^\infty}(k_c)
+k_c\partial_kS_{\rho_c^\infty}(k_c)
\right]
&=2\alpha k_c.
\end{aligned}
\end{equation}

A solution of \eqref{eq:wl-critical-pair} is only a candidate
for the first loss of linear stability. It must correspond to
a global maximum, so
\[
\sigma(k;\rho_c^\infty,\lambda)\le0
\quad\text{for every }k\in\mathbb R,
\]
with equality at $k_c$. The radius must also be the one specified by \eqref{eq:wl-critical-radius}. No monotonicity in $\rho$ is assumed, so instability need not
persist at every radius above $\rho_c^\infty$.

The wavelength associated with the critical mode is
\begin{equation*}
\ell_c=\frac{2\pi}{k_c}.
\end{equation*}
In an unstable regime, a fastest growing wavenumber satisfies
\[
k_*(\rho,\lambda)
\in
\operatorname*{arg\,max}_{k\ge0}\sigma(k;\rho,\lambda),
\]
and its wavelength is $\ell_*=2\pi/k_*$. More than one maximizing
wavenumber may exist. These wavelengths describe the spatial
scales predicted during the initial growth of small
perturbations; they do not determine the final nonlinear
cell density profile.

For the kernel used in this work,
\begin{equation}
\label{eq:wl-arctan-kernel}
h(\xi)
=
\frac{0.1\arctan(0.2\xi)}{\arctan(2)},
\quad |\xi|\le\rho,
\end{equation}
the denominator remains fixed as $\rho$ varies. Define
\[
A_\rho(k)
=
\int_0^\rho
\arctan(0.2\xi)\sin(k\xi)\,d\xi,
\quad
B_\rho(k)
=
\int_0^\rho
\xi\arctan(0.2\xi)\cos(k\xi)\,d\xi.
\]
Then
\[
S_\rho(k)=\frac{0.2}{\arctan(2)}A_\rho(k),
\quad
\partial_kS_\rho(k)=\frac{0.2}{\arctan(2)}B_\rho(k),
\]
and \eqref{eq:wl-critical-pair} becomes
\begin{equation*}
\begin{aligned}
\frac{0.2(\lambda-2)}{\arctan(2)}
\,k_c A_{\rho_c^\infty}(k_c)
&=\alpha k_c^2+1,\\
\frac{0.2(\lambda-2)}{\arctan(2)}
\left[
A_{\rho_c^\infty}(k_c)
+k_cB_{\rho_c^\infty}(k_c)
\right]
&=2\alpha k_c.
\end{aligned}
\end{equation*}
These equations determine the critical pair implicitly,
subject to the global maximum and first transition checks
above. The integrals and the critical pair are evaluated
numerically.

The whole line results provide a reference for the bounded
interval. They do not give the exact threshold for the
prescribed healthy exterior, which requires analysis of the
bounded linearized operator.

\section{Bounded domain formulation and stability}
\label{sec:bounded-formulation}

We now consider equation~\eqref{eq:wound_healing_model} on the
bounded interval $\Omega=(-L,L)$. For a sensing radius $\rho>0$, define the exterior collar by
\begin{equation*}
\Omega_\rho=(-L-\rho,-L]\cup[L,L+\rho).
\end{equation*}
The cell density is maintained at the healed value $1$ on this
collar. Thus, cells inside $\Omega$ can interact with surrounding
cells whose density is prescribed. To distinguish the unknown interior density from its prescribed
exterior values, define
\begin{equation}
\label{eq:bd-healthy-extension}
n^e(x,t)=
\begin{cases}
n(x,t), & x\in(-L,L),\\
1, & x\notin(-L,L).
\end{cases}
\end{equation}
Only the values of this extension in $\Omega\cup\Omega_\rho$
are sampled by the adhesion integral. Extending the density by
$1$ beyond the collar is therefore a convenient convention that
does not change the equation on $\Omega$. For $x\in\Omega$, define the adhesion velocity by
\begin{equation*}
V[n](x,t)
=\int_{-\rho}^{\rho}g\bigl(n^e(x+\xi,t)\bigr)h(\xi)\,d\xi.
\end{equation*}
The integration interval remains $[-\rho,\rho]$ for every
$x\in\Omega$. When $x+\xi$ lies outside $\Omega$, the integrand
uses the prescribed healthy density.

With $\alpha>0$ and $f(n)=n(1-n)$, the bounded problem is
\begin{equation}
\label{eq:bd-problem}
\begin{cases}
n_t=\alpha n_{xx}-\partial_x\bigl(nV[n]\bigr)+n(1-n), &x\in(-L,L),\quad t>0,\\
n(-L,t)=n(L,t)=1,  &t\ge0,\\
n(x,0)=n_0(x), &x\in(-L,L),
\end{cases}
\end{equation}
We assume that the initial density is nonnegative and satisfies the condition
$n_0(\pm L)=1.$ The endpoint conditions in \eqref{eq:bd-problem} specify the Dirichlet boundary values for the local diffusion term. They are understood as interior traces and agree with the healthy exterior prescription.

The parameter $L$ is the domain half-length; it is not necessarily the half-width of the initial wound. Keeping the initial wound width fixed while varying $L$ changes the distance between the wound and the prescribed healthy exterior.

For the estimates in this section, assume that $h$ is real and odd and that
$h\in W^{1,1}(-\rho,\rho)$. This additional regularity permits integration by
parts in the sensing variable.

The following calculations concern sufficiently regular
solutions on their interval of existence. In particular, on
each compact time interval $[0,T]$ contained in that interval,
we assume
\[
n-1\in C\bigl([0,T];H_0^1(\Omega)\bigr),
\]
with sufficient time regularity to justify the energy identity
below.

At $n\equiv1$, 
\[
V[1]=g(1)\int_{-\rho}^{\rho}h(\xi)\,d\xi=0.
\]
The diffusion and logistic terms also vanish. Thus, $n\equiv1$ is an equilibrium solution of \eqref{eq:bd-problem} that satisfies both the boundary conditions and the prescribed exterior density.

\subsection{A bound on the velocity derivative}
\label{subsec:bd-velocity-bound}

To control the cell density, we first estimate
$\partial_xV[n]$. For sufficiently regular $n$, integration
by parts in $\xi$ gives
\begin{align}
\partial_xV[n](x,t)
&=\int_{-\rho}^{\rho}
\partial_\xi\!\left[g\bigl(n^e(x+\xi,t)\bigr)\right]
h(\xi)\,d\xi \notag\\
&=\left[
g\bigl(n^e(x+\xi,t)\bigr)h(\xi)
\right]_{-\rho}^{\rho}
-\int_{-\rho}^{\rho}
g\bigl(n^e(x+\xi,t)\bigr)h'(\xi)\,d\xi \notag\\
&=h(\rho)\left[
g\bigl(n^e(x+\rho,t)\bigr)
+g\bigl(n^e(x-\rho,t)\bigr)
\right] \notag\\
&\hspace{1em}
-\int_{-\rho}^{\rho}
g\bigl(n^e(x+\xi,t)\bigr)h'(\xi)\,d\xi.
\label{eq:bd-velocity-derivative}
\end{align}
The extension $n^e$ is continuous at $\pm L$ because the interior traces equal $1$. Thus there are no jump terms at the endpoints of $\Omega$ in this calculation. 
Define
\begin{equation*}
J_\rho
=2|h(\rho)|+\int_{-\rho}^{\rho}|h'(\xi)|\,d\xi.
\end{equation*}
If $n^e$ takes values in an interval $I$ containing $1$, then
\eqref{eq:bd-velocity-derivative} implies
\begin{equation}
\label{eq:bd-general-velocity-bound}
\|\partial_xV[n](\cdot,t)\|_{L^\infty(\Omega)}
\le J_\rho\max_{s\in I}|g(s)|.
\end{equation}

For the arctangent kernel \eqref{eq:wl-arctan-kernel} used in this work,
\[
h(\xi)=\frac{0.1\arctan(0.2\xi)}{\arctan(2)},
\quad |\xi|\le\rho,
\]
we have
\[
h'(\xi)=\frac{0.02}{\arctan(2)(1+0.04\xi^2)}>0,
\]
and therefore
\begin{equation*}
J_\rho=4h(\rho)
=\frac{0.4\arctan(0.2\rho)}{\arctan(2)}.
\end{equation*}
In particular, $h(\rho)\ne0$ for $\rho>0$, so omitting the
endpoint terms would give an incorrect derivative formula.

\subsection{An invariant interval for the cell density}
\label{subsec:bd-invariant-interval}

The following lemma gives a sufficient condition under which
the density remains in a prescribed positive bounded interval.

\begin{lemma}[Invariant density interval]
\label{lem:bd-invariant-interval}
Let $0<\delta_0<1<M$, and set
\begin{equation*}
G_{\max}=\max_{s\in[\delta_0,M]}|g(s)|,
\quad K_\rho=G_{\max}J_\rho.
\end{equation*}
Suppose that
\begin{equation*}
\delta_0\le n_0(x)\le M
\quad\text{for every }x\in[-L,L],
\end{equation*}
and
\begin{equation}
\label{eq:bd-invariant-condition}
K_\rho\le\min\{1-\delta_0,M-1\}.
\end{equation}
Then every sufficiently regular solution of
\eqref{eq:bd-problem} satisfies
\begin{equation*}
\delta_0\le n(x,t)\le M
\end{equation*}
throughout its interval of existence.
\end{lemma}

\begin{proof}
Expanding the flux derivative gives
\begin{equation}
\label{eq:bd-expanded-density-equation}
n_t=\alpha n_{xx}-V[n]n_x
+n\bigl(1-n-\partial_xV[n]\bigr).
\end{equation}
Fix a finite time $T$
strictly below the endpoint of the solution's existence
interval. Define
\[
\widehat L_g
=\max_{s\in[\delta_0/2,M+1]}|g'(s)|,
\quad D=J_\rho\widehat L_g,
\quad C=1+(M+1)D.
\]
Choose $\varepsilon>0$ sufficiently small that
\[
q(t):=\varepsilon e^{Ct}
<\min\{\delta_0/2,1\}
\quad\text{for }0\le t\le T,
\]
and introduce the moving barriers
\[
a(t)=\delta_0-q(t),
\quad b(t)=M+q(t).
\]
Initially, $a(0)<n_0<b(0)$. The prescribed boundary and
exterior value $1$ also lies strictly between these barriers.

Suppose that a first contact with either barrier occurs at
an interior point $(x_*,t_*)$, where $0<t_*\le T$. Up to
this time, all sampled density values lie in
$[\delta_0-q(t),M+q(t)]$. Every point in this interval is
within distance $q(t)$ of $[\delta_0,M]$. The mean value
theorem therefore gives
\[
\max_{s\in[\delta_0-q(t),M+q(t)]}|g(s)|
\le G_{\max}+\widehat L_gq(t),
\]
and hence, by \eqref{eq:bd-general-velocity-bound}
\begin{equation}
\label{eq:bd-expanded-band-velocity-bound}
|\partial_xV[n](x,t)|\le K_\rho+Dq(t)
\end{equation}
up to the first contact.

At a lower contact, $n=a(t_*)$, $n_x=0$, $n_{xx}\ge0$,
and $n_t\le a'(t_*)=-Cq(t_*)$. On the other hand,
\eqref{eq:bd-expanded-density-equation} and
\eqref{eq:bd-expanded-band-velocity-bound} imply, with
$q=q(t_*)$,
\[
\begin{aligned}
n_t
&\ge a(t_*)\bigl[1-\delta_0+q-K_\rho-Dq\bigr]\\
&\ge-a(t_*)Dq\\
&\ge-(M+1)Dq
>-Cq,
\end{aligned}
\]
a contradiction.

At an upper contact, $n=b(t_*)$, $n_x=0$, $n_{xx}\le0$,
and $n_t\ge b'(t_*)=Cq(t_*)$. However, \eqref{eq:bd-expanded-density-equation} and
\eqref{eq:bd-expanded-band-velocity-bound} imply, with
$q=q(t_*)$,
\[
\begin{aligned}
n_t
&\le b(t_*)\bigl[1-M-q+K_\rho+Dq\bigr]\\
&\le b(t_*)Dq\\
&\le(M+1)Dq
<Cq,
\end{aligned}
\]
which is again a contradiction.

Thus $a(t)<n(x,t)<b(t)$ on $[0,T]$. Letting
$\varepsilon\downarrow0$ proves
$\delta_0\le n(x,t)\le M$ there. Since $T$ was arbitrary,
the bounds hold throughout the solution's existence interval.
\end{proof}

Condition~\eqref{eq:bd-invariant-condition} is sufficient and
can be restrictive. Its failure does not imply that the
density leaves the interval. The positive lower bound also
excludes an initially empty part of the wound; such initial
data require a separate positivity argument before an estimate
using $\delta_0>0$ can be applied.

\subsection{Nonlinear exponential convergence}
\label{subsec:bd-nonlinear-stability}

For the nonlinear estimates, let $0<\delta_0<1<M$ and assume that \eqref{eq:bd-problem}
admits a global sufficiently regular solution satisfying
\begin{equation}
\label{eq:bd-density-bounds}
\delta_0\le n(x,t)\le M \quad (x\in\Omega,\ t\ge0).
\end{equation}
By Lemma~\ref{lem:bd-invariant-interval}, \eqref{eq:bd-density-bounds} holds
whenever $\delta_0\le n_0\le M$ and \eqref{eq:bd-invariant-condition} is
satisfied; outside that parameter range it is an assumption.

Write $\|\cdot\|_2=\|\cdot\|_{L^2(\Omega)}$ and recall
\[
H_\rho=\int_{-\rho}^{\rho}|h(\xi)|\,d\xi.
\]
For a fixed density interval $[\delta_0,M]$, define
\begin{equation*}
L_g=\max_{s\in[\delta_0,M]}|g'(s)|
=\max\{|\lambda-2\delta_0|,|\lambda-2M|\},
\quad C_\rho=ML_gH_\rho.
\end{equation*}
Here $L_g$ is a Lipschitz constant for $g$ and is distinct
from the domain half-length $L$.

Let $w=n-1$, $w_0=n_0-1$, and let $w^e$ be the zero
extension of $w$ outside $\Omega$. Then $n^e=1+w^e$ and
$w(\pm L,t)=0$. Since $h$ is odd,
\begin{equation*}
V[n](x,t)=\int_{-\rho}^{\rho}
\bigl[g(1+w^e(x+\xi,t))-g(1)\bigr]h(\xi)\,d\xi.
\end{equation*}
Whenever $\delta_0\le n\le M$, the mean value theorem gives
\[
|g(1+w^e)-g(1)|\le L_g|w^e|.
\]
Young's convolution inequality, followed by restriction to
$\Omega$, therefore yields
\begin{equation}
\label{eq:bd-velocity-l2-bound}
\|V[n](\cdot,t)\|_2
\le L_gH_\rho\|w(\cdot,t)\|_2.
\end{equation}

The exact perturbation equation is
\begin{equation*}
w_t=\alpha w_{xx}-\partial_x(nV[n])-nw,
\quad w(\pm L,t)=0.
\end{equation*}
Multiplication by $w$ and integration by parts give
\begin{equation*}
\frac12\frac{d}{dt}\|w\|_2^2
=-\alpha\|w_x\|_2^2
+\int_\Omega nV[n]w_x\,dx
-\int_\Omega nw^2\,dx.
\end{equation*}
There are no boundary terms because $w$ has zero trace.
Using \eqref{eq:bd-velocity-l2-bound} and the density bounds,
we obtain
\begin{equation}
\label{eq:bd-energy-inequality}
\frac12\frac{d}{dt}\|w\|_2^2
\le-\alpha\|w_x\|_2^2
+C_\rho\|w_x\|_2\|w\|_2-\delta_0\|w\|_2^2.
\end{equation}

\begin{theorem}[Nonlinear exponential convergence]
\label{thm:bd-nonlinear-decay}
Assume that a global sufficiently regular solution of
\eqref{eq:bd-problem} satisfies
\[
0<\delta_0\le n(x,t)\le M
\quad (x\in\Omega,\ t\ge0).
\]
If
\begin{equation}
\label{eq:bd-nonlinear-condition}
C_\rho^2<4\alpha\delta_0,
\end{equation}
then
\begin{equation}
\label{eq:bd-nonlinear-decay}
\|n(\cdot,t)-1\|_2
\le e^{-\mu t}\|n_0-1\|_2,
\quad
\mu=\delta_0-\frac{C_\rho^2}{4\alpha}>0.
\end{equation}
\end{theorem}

\begin{proof}
For $X=\|w_x\|_2$ and $Y=\|w\|_2$, completing the square gives
\[
-\alpha X^2+C_\rho XY-\delta_0Y^2
=-\alpha\left(X-\frac{C_\rho}{2\alpha}Y\right)^2
-\left(\delta_0-\frac{C_\rho^2}{4\alpha}\right)Y^2.
\]
It follows from \eqref{eq:bd-energy-inequality} that
\[
\frac{d}{dt}\|w\|_2^2\le-2\mu\|w\|_2^2.
\]
Gronwall's inequality and taking square roots prove
\eqref{eq:bd-nonlinear-decay}.
\end{proof}

The invariant interval in Lemma~\ref{lem:bd-invariant-interval}
provides one sufficient way to obtain the density bounds
during a solution's existence. Applying
Theorem~\ref{thm:bd-nonlinear-decay} for all $t\ge0$ also
requires global existence. The theorem concerns the full
nonlinear equation and does not require the initial
perturbation to be infinitesimal. Its hypotheses do require
a uniform positive lower bound on the density.

\begin{remark}
Condition \eqref{eq:bd-invariant-condition} restricts the sensing
radius. Since $\delta_0<1<M$, the value $s=1$ lies in $[\delta_0,M]$,
so $G_{\max}\ge g(1)=\lambda-1$, and \eqref{eq:bd-invariant-condition}
requires $(\lambda-1)J_\rho<1$. For the kernel \eqref{eq:wl-arctan-kernel},
$J_\rho=4h(\rho)$ increases in $\rho$ with
$\sup_{\rho>0}J_\rho=0.2\pi/\arctan 2\approx0.568$, so
\eqref{eq:bd-invariant-condition} is satisfiable for every $\rho>0$
when $\lambda\lesssim2.6$, and for $\rho\lesssim2.5$ when $\lambda=5$.
With $\alpha=1$, $\lambda=5$, $\rho=1$, $\delta_0=0.25$ and $M=1.5$,
both \eqref{eq:bd-invariant-condition} and
\eqref{eq:bd-nonlinear-condition} hold, giving unconditional
exponential convergence at rate $\mu\approx0.246$. The parameter
values used in Section~\ref{sec:numerical-experiments} lie outside this range, so the density
bounds are assumed rather than derived.
\end{remark}

\subsection{Domain dependence and the sufficient length bound}
\label{subsec:bd-length-bound}

The preceding estimate does not use the domain size.
To retain its effect, set
\begin{equation*}
\kappa_L=\frac{\pi}{2L},
\quad
\|v_x\|_2^2\ge\kappa_L^2\|v\|_2^2
\quad\text{for }v\in H_0^1(-L,L).
\end{equation*}
Define
\begin{equation*}
\mu_L=
\begin{cases}
\displaystyle
\delta_0-\frac{C_\rho^2}{4\alpha},
&\displaystyle \kappa_L\le\frac{C_\rho}{2\alpha},\\[2mm]
\displaystyle
\alpha\kappa_L^2-C_\rho\kappa_L+\delta_0,
&\displaystyle \kappa_L>\frac{C_\rho}{2\alpha}.
\end{cases}
\end{equation*}

\begin{theorem}[Domain dependent nonlinear decay]
\label{thm:bd-domain-decay}
Assume that a global sufficiently regular solution remains
in $[\delta_0,M]$, where $0<\delta_0<1<M$.
If $\mu_L>0$, then
\begin{equation*}
\|n(\cdot,t)-1\|_2
\le e^{-\mu_Lt}\|n_0-1\|_2,
\quad t\ge0.
\end{equation*}
This conclusion does not require
\eqref{eq:bd-nonlinear-condition}.
\end{theorem}

\begin{proof}
For nonzero $w\in H_0^1(\Omega)$, put
\[
r=\frac{\|w_x\|_2}{\|w\|_2}\ge\kappa_L.
\]
Equation~\eqref{eq:bd-energy-inequality} becomes
\[
\frac12\frac{d}{dt}\|w\|_2^2
\le\bigl(-\alpha r^2+C_\rho r-\delta_0\bigr)\|w\|_2^2.
\]
The concave quadratic
\[
q(r)=-\alpha r^2+C_\rho r-\delta_0
\]
has its maximum at $r_*=C_\rho/(2\alpha)$.
Consequently,
\[
\sup_{r\ge\kappa_L}q(r)
=
\begin{cases}
\displaystyle \frac{C_\rho^2}{4\alpha}-\delta_0,
&\kappa_L\le r_*,\\[1mm]
-\alpha\kappa_L^2+C_\rho\kappa_L-\delta_0,
&\kappa_L>r_*,
\end{cases}
=-\mu_L.
\]
The resulting differential inequality also holds when $w=0$,
since \eqref{eq:bd-energy-inequality} then has zero right-hand
side. Gronwall's inequality proves the assertion.
\end{proof}

\begin{corollary}[Sufficient domain half-length for healing]
\label{cor:bd-Lstab}
Fix $\alpha$, $\lambda$, $\rho$, $\delta_0$, and $M$, and define
\begin{equation}
\label{eq:bd-Lstab}
L_{\mathrm{stab}}(\lambda,\rho)=
\begin{cases}
+\infty,
& C_\rho^2<4\alpha\delta_0,\\[2mm]
\displaystyle
\frac{\pi\alpha}
{C_\rho+\sqrt{C_\rho^2-4\alpha\delta_0}},
& C_\rho^2\ge4\alpha\delta_0.
\end{cases}
\end{equation}
Under the global solution and density bound assumptions of
Theorem~\ref{thm:bd-domain-decay}, the condition
\[
L<L_{\mathrm{stab}}(\lambda,\rho)
\]
implies exponential convergence to $n\equiv1$ with rate
$\mu_L>0$.
\end{corollary}

\begin{proof}
If $C_\rho^2<4\alpha\delta_0$, then $q(r)<0$ for every
$r\ge0$, so $\mu_L>0$ for every $L>0$.

If $C_\rho^2\ge4\alpha\delta_0$, the larger root of $q$ is
\[
r_+=\frac{C_\rho+\sqrt{C_\rho^2-4\alpha\delta_0}}{2\alpha}.
\]
Since $r_+\ge r_*$, the supremum of $q$ over
$[\kappa_L,\infty)$ is strictly negative precisely when
$\kappa_L>r_+$. Substituting $\kappa_L=\pi/(2L)$ gives the
finite bound in \eqref{eq:bd-Lstab}.
\end{proof}

The dependence on $\alpha$, $\delta_0$, and $M$ is suppressed
in the notation $L_{\mathrm{stab}}(\lambda,\rho)$.
The bound concerns the domain half-length; the corresponding
full interval length is $2L_{\mathrm{stab}}$.
Failure of this sufficient test does not establish failure
to heal. When the finite bound is attained, this estimate
gives $\mu_L=0$ and does not establish exponential decay.

For the arctangent kernel \ref{eq:wl-arctan-kernel},
\begin{equation*}
H_\rho=
\frac{0.2}{\arctan(2)}
\left[
\rho\arctan(0.2\rho)
-\frac52\log(1+0.04\rho^2)
\right].
\end{equation*}
Thus all quantities in \eqref{eq:bd-Lstab} can be evaluated
explicitly once the density interval and model parameters
have been specified.

\subsection{The bounded linearized operator}
\label{subsec:bd-linear-spectrum}

For $v\in L^2(\Omega)$, let $v^e$ be its zero extension and set
\begin{equation*}
(\mathcal V_{L,\rho}v)(x)
=\int_{-\rho}^{\rho}v^e(x+\xi)h(\xi)\,d\xi.
\end{equation*}
Since
\[
g(1+w^e)-g(1)=(\lambda-2)w^e-(w^e)^2,
\]
retaining the terms linear in $w$ gives
\begin{equation*}
w_t=\mathcal A_{L,\rho,\lambda}w,
\quad w(\pm L,t)=0,
\end{equation*}
where
\begin{equation}
\label{eq:bd-linearized-operator}
\begin{aligned}
\mathcal A_{L,\rho,\lambda}v
&=\alpha v''-v+(\lambda-2)\mathcal K_{L,\rho}v,\\
\mathcal K_{L,\rho}v
&=-\partial_x(\mathcal V_{L,\rho}v),\\
D(\mathcal A_{L,\rho,\lambda})
&=H^2(\Omega)\cap H_0^1(\Omega).
\end{aligned}
\end{equation}
Integration by parts in the sensing variable yields
\begin{equation}
\label{eq:bd-K-explicit}
(\mathcal K_{L,\rho}v)(x)
=\int_{-\rho}^{\rho}h'(\xi)v^e(x+\xi)\,d\xi
-h(\rho)\bigl[v^e(x+\rho)+v^e(x-\rho)\bigr].
\end{equation}
This formula initially holds for $v\in H_0^1(\Omega)$ and
defines a bounded operator on $L^2(\Omega)$, with
\begin{equation}
\label{eq:bd-K-norm}
\|\mathcal K_{L,\rho}v\|_2\le J_\rho\|v\|_2.
\end{equation}

\begin{proposition}[Bounded spectral characterization]
\label{prop:bd-spectrum}
The operator $\mathcal A_{L,\rho,\lambda}$ is self-adjoint
on $L^2(\Omega)$ and has compact resolvent. Its eigenvalues,
counted with multiplicity, are real and satisfy
\[
\sigma_1(L;\rho,\lambda)\ge\sigma_2(L;\rho,\lambda)
\ge\cdots\longrightarrow-\infty.
\]
Define the largest bounded growth rate by
\begin{equation*}
s_L(\rho;\lambda)
=\max\operatorname{spec}(\mathcal A_{L,\rho,\lambda})
=\sigma_1(L;\rho,\lambda).
\end{equation*}
For real-valued test functions,
\begin{equation}
\label{eq:bd-rayleigh}
s_L(\rho;\lambda)
=\sup_{\substack{v\in H_0^1(\Omega)\\v\ne0}}
\frac{
-\alpha\|v'\|_2^2-\|v\|_2^2
+(\lambda-2)\displaystyle\int_\Omega
v'(x)(\mathcal V_{L,\rho}v)(x)\,dx
}{\|v\|_2^2}.
\end{equation}
The linearized solution satisfies
\begin{equation}
\label{eq:bd-linear-growth}
\|w(\cdot,t)\|_2
\le e^{s_L(\rho;\lambda)t}\|w_0\|_2.
\end{equation}
\end{proposition}

\begin{proof}
Since $h$ is odd, $h'$ is even almost everywhere.
After writing $y=x+\xi$, the integral term in
\eqref{eq:bd-K-explicit} has the symmetric kernel
\[
h'(y-x)\,\mathbf 1_{\{|y-x|<\rho\}}.
\]
Moreover, the zero-extended shifts
\[
T_+v(x)=v^e(x+\rho),
\quad T_-v(x)=v^e(x-\rho)
\]
satisfy $T_+^*=T_-$. Thus $\mathcal K_{L,\rho}$ is bounded
and self-adjoint.

The Dirichlet operator $\alpha\partial_{xx}-I$ is
self-adjoint with compact resolvent. Its bounded
self-adjoint perturbation
$(\lambda-2)\mathcal K_{L,\rho}$ preserves self-adjointness
on the same domain. The resolvent identity preserves
compactness of the resolvent. The bound
\eqref{eq:bd-K-norm} and the Dirichlet eigenvalues also show
that the perturbed eigenvalues tend to $-\infty$.

For $v\in D(\mathcal A_{L,\rho,\lambda})$,
integration by parts gives
\[
\langle v,\mathcal A_{L,\rho,\lambda}v\rangle
=-\alpha\|v'\|_2^2-\|v\|_2^2
+(\lambda-2)\int_\Omega v'\mathcal V_{L,\rho}v\,dx.
\]
The quadratic form extends continuously to $H_0^1(\Omega)$.
The variational characterization of the largest eigenvalue
gives \eqref{eq:bd-rayleigh}.

Let $\{\phi_j\}_{j\ge1}$ be an orthonormal eigenbasis.
For $w_0=\sum_jc_j\phi_j$, the linearized solution is
\[
w(\cdot,t)=\sum_{j\ge1}c_je^{\sigma_jt}\phi_j.
\]
Hence
\[
\|w(\cdot,t)\|_2^2
=\sum_{j\ge1}|c_j|^2e^{2\sigma_jt}
\le e^{2s_Lt}\|w_0\|_2^2,
\]
which proves \eqref{eq:bd-linear-growth}.
\end{proof}

For fixed $L$, $\lambda$, and $\alpha$, we write
\begin{equation}\label{eq:bd-s-rho}
s(\rho)=s_L(\rho;\lambda).
\end{equation}
If $s(\rho)<0$, every linearized perturbation decays
exponentially. If $s(\rho)>0$, an eigenfunction associated
with $s(\rho)$ gives an exponentially growing perturbation.
If $s(\rho)=0$, the linearized equation has a nonzero
stationary eigenfunction; its solution operator is bounded,
but not every perturbation decays to zero. These spectral statements alone do not determine whether
solutions of~\eqref{eq:bd-problem} converge to the healed
equilibrium $n\equiv1$.

\subsection{Critical sensing radius and linear length threshold}
\label{subsec:bd-linear-thresholds}

For fixed $L$, $\alpha$, $\lambda$, and a specified kernel
family, define
\begin{equation}
\label{eq:bd-critical-radius}
\rho_c^{(L)}(\lambda)
=\inf\{\rho>0:s_L(\rho;\lambda)>0\}.
\end{equation}
The infimum of the empty set is taken to be $+\infty$.
All dependence of the kernel on $\rho$, including its
normalization, is included in the definition of the operator.

\begin{proposition}[Growth rate at the critical radius]
\label{prop:bd-critical-growth}
Suppose that $0<\rho_c^{(L)}(\lambda)<\infty$ and that
$s_L(\cdot;\lambda)$ is continuous at this radius. Then
\begin{equation}
\label{eq:bd-critical-root}
s_L\bigl(\rho_c^{(L)}(\lambda);\lambda\bigr)=0.
\end{equation}
Consequently, there is a nonzero critical eigenfunction
$\phi_c\in H^2(\Omega)\cap H_0^1(\Omega)$ satisfying
\begin{equation*}
\mathcal A_{L,\rho_c^{(L)}(\lambda),\lambda}\phi_c=0.
\end{equation*}
\end{proposition}

\begin{proof}
If the value of $s_L$ at the critical radius were positive,
continuity would give a smaller radius with positive $s_L$.
If that value were negative, continuity would exclude
positive values in a neighborhood of the infimum, although
such values must occur arbitrarily close to it.
Both alternatives contradict \eqref{eq:bd-critical-radius}.
The eigenfunction follows from
Proposition~\ref{prop:bd-spectrum}.
\end{proof}

A sufficient condition for the required continuity is
$L^1(\mathbb R)$-continuity of the zero-extended kernel:
\begin{equation}
\label{eq:bd-kernel-family-continuity}
\|\widetilde h_\rho-\widetilde h_{\rho_0}\|_{L^1(\mathbb R)}
\longrightarrow0
\quad\text{as }\rho\longrightarrow\rho_0,
\quad
\widetilde h_\rho(\xi)
=h_\rho(\xi)\mathbf 1_{\{|\xi|<\rho\}}.
\end{equation}
Here $h_\rho$ makes explicit any dependence of the kernel
on its sensing radius.

To justify this assertion, for $\|v\|_2=1$ the difference
between the two Rayleigh forms is bounded by
\[
|\lambda-2|\,
\|\widetilde h_\rho-\widetilde h_{\rho_0}\|_{L^1(\mathbb R)}
\|v'\|_2.
\]
The derivative norms of normalized maximizing eigenfunctions
are locally bounded in $\rho$: testing the Rayleigh quotient
with one fixed function gives a local lower bound, while
\[
s_L(\rho;\lambda)
\le-\alpha r^2-1+|\lambda-2|H_\rho r,
\quad r=\|v'\|_2,
\]
at a normalized maximizer bounds $r$ from above.
Comparing the two maximizing eigenfunctions in the
variational formula proves continuity.
Condition~\eqref{eq:bd-kernel-family-continuity} holds for
the kernel \eqref{eq:wl-arctan-kernel}.

Equation~\eqref{eq:bd-critical-root} is the scalar equation
used to locate a candidate threshold. Definition
\eqref{eq:bd-critical-radius} specifies which root to select.
No monotonicity of $s_L$ in $\rho$ is assumed; a zero that
merely touches the axis need not mark the onset of positive
growth rates.

For the fixed shape arctangent kernel, the perturbation
operator is independent of $\rho$ once $\rho\ge2L$.
Indeed,
\[
(\mathcal V_{L,\rho}v)(x)
=\int_{-L}^{L}
v(y)h(y-x)\mathbf 1_{\{|y-x|<\rho\}}\,dy,
\]
and $|y-x|<2L$ for $x,y\in\Omega$.
This observation does not apply unchanged to families whose
shape or normalization varies with $\rho$.

\begin{proposition}[Dependence on domain size]
\label{prop:bd-domain-monotonicity}
Fix $\alpha$, $\rho$, $\lambda$, and a kernel that does not
depend on $L$. Then $L\mapsto s_L(\rho;\lambda)$ is
nondecreasing. Moreover,
\begin{equation}
\label{eq:bd-small-domain-spectral-bound}
s_L(\rho;\lambda)
\le-\alpha\left(\frac{\pi}{2L}\right)^2
-1+|\lambda-2|J_\rho,
\end{equation}
so sufficiently small domains are linearly exponentially
stable.
\end{proposition}

\begin{proof}
Let $0<L_1<L_2$. Zero extension maps
$H_0^1(-L_1,L_1)$ into $H_0^1(-L_2,L_2)$
and preserves the Rayleigh quotient.
Taking suprema in \eqref{eq:bd-rayleigh} proves
$s_{L_1}\le s_{L_2}$.
The operator bound \eqref{eq:bd-K-norm} and the Poincar\'e
inequality give \eqref{eq:bd-small-domain-spectral-bound},
whose right-hand side tends to $-\infty$ as $L\downarrow0$.
\end{proof}

To describe loss of strict exponential decay as $L$ increases, define
\begin{equation}
\label{eq:bd-Llin}
L_c(\rho,\lambda)
=\inf\{L>0:s_L(\rho;\lambda)\ge0\},
\end{equation}
again taking the infimum of the empty set to be $+\infty$.
For $L<L_c$, the linearized equation is
exponentially stable. If $L_c<\infty$, then
$s_L\ge0$ for $L>L_c$, by monotonicity.
Positive growth still requires $s_L>0$; the definition
includes neutral cases with $s_L=0$. If the threshold is finite
and $L\mapsto s_L(\rho;\lambda)$ is continuous there, then
\[
s_{L_c}(\rho;\lambda)=0.
\]

The three quantities have different roles:
$L_{\mathrm{stab}}$ is a sufficient bound for exponential
healing under the stated density and global existence
assumptions; $L_c$ marks the loss of strict linear exponential
stability as the domain half-length increases; and
$\rho_c^{(L)}$ marks the onset of a positive linear
growth rate as the sensing radius varies.

\subsection{Critical eigenfunctions and whole line comparison}
\label{subsec:bd-critical-modes}

The normalized Dirichlet sine basis is
\[
e_j(x)=\frac{1}{\sqrt L}
\sin\left(\frac{j\pi(x+L)}{2L}\right),
\quad
k_j=\frac{j\pi}{2L},
\quad j=1,2,\ldots.
\]
Although these functions diagonalize the local diffusion
operator, they generally do not diagonalize
$\mathcal K_{L,\rho}$. In this basis, the coupled
Galerkin matrix has entries
\begin{equation}
\label{eq:bd-galerkin-matrix}
A_{ij}
=-(\alpha k_j^2+1)\delta_{ij}
+(\lambda-2)\int_\Omega
e_i'(x)(\mathcal V_{L,\rho}e_j)(x)\,dx.
\end{equation}
The matrix is symmetric when the integrals are evaluated
exactly. The largest eigenvalue of its truncation supplies
an approximation of $s_L(\rho;\lambda)$.

At a bounded threshold, first compute $\phi_c$ and then
expand it as
\[
\phi_c=\sum_{j\ge1}a_je_j,
\quad
a_j=\int_\Omega\phi_c(x)e_j(x)\,dx.
\]
If one component dominates, its spatial scale may be
summarized by
\begin{equation}
\label{eq:bd-dominant-wavenumber}
j_{\mathrm{dom}}\in\operatorname*{arg\,max}_{j\ge1}|a_j|^2,
\quad
k_{c,\mathrm{dom}}^{(L)}=\frac{j_{\mathrm{dom}}\pi}{2L},
\quad
\ell_{c,\mathrm{dom}}^{(L)}
=\frac{2\pi}{k_{c,\mathrm{dom}}^{(L)}}.
\end{equation}
This is a description of the computed eigenfunction, not
an assumption that it is a single sine mode.

Reflection symmetry leaves the even and odd subspaces
invariant. The largest eigenvalue must therefore be sought
over both subspaces unless perturbations are explicitly
restricted to one symmetry class.

For comparison with the whole line problem, the
zero extended Rayleigh quotient can be written using
Plancherel's identity as
\begin{equation}
\label{eq:bd-whole line-rayleigh}
\mathcal R_L(v)
=
\frac{\displaystyle\int_{\mathbb R}
\sigma(k;\rho,\lambda)|\widehat{v^e}(k)|^2\,dk}
{\displaystyle\int_{\mathbb R}
|\widehat{v^e}(k)|^2\,dk}.
\end{equation}
Here $\sigma$ is the whole line dispersion relation for the
same kernel and parameters. Hence
\begin{equation*}
s_L(\rho;\lambda)\le F(\rho;\lambda).
\end{equation*}
In particular, whole line linear exponential stability
implies bounded linear exponential stability for this
healthy exterior formulation. The converse need not hold.

Comparing $k_j$ with the whole line unstable band gives
a qualitative indication of how domain size affects stability.
The bounded thresholds, however, are defined by the full
operator and its largest eigenvalue. Numerical evaluation
of these thresholds and comparisons with finite amplitude
wound simulations are addressed in the numerical section.

\section{Numerical experiments}
\label{sec:numerical-experiments}

We compare the spectral stability predictions with recovery
from finite amplitude wounds. The computations determine the
critical sensing radii $\rho_c^\infty$ and $\rho_c^{(L)}$, the
critical domain half-length $L_c$, and the dependence of healing
on domain size. Throughout, $\alpha=1$ and the arctangent kernel
\eqref{eq:wl-arctan-kernel} retains its fixed normalization.

\subsection{Numerical method}
\label{subsec:num-method}

To solve~\eqref{eq:bd-problem}, we discretize the
domain using a uniform grid with step size $\Delta x$ and approximate the diffusion term by Fourier continuation~\cite{albin2014discrete}. The prescribed exterior collar values
are used in the adhesion integral and are kept separate
from the continuation used for differentiation.

At each physical grid point $x_i$, we approximate the adhesion
velocity by piecewise-linear product integration:
\[
V[n](x_i,t)\approx V_i(t)
:=\sum_{j=-m}^{m}\omega_j
\bigl[g(n^e_{i+j}(t))-g(1)\bigr],
\quad
m=\left\lceil\frac{\rho}{\Delta x}\right\rceil.
\]
Here $n^e_{i+j}(t)$ is the density at $x_i+j\Delta x$,
including the prescribed collar values. The weights are
defined by
\[
\omega_j=\int_{-\rho}^{\rho}B_j(\xi)h(\xi)\,d\xi,
\quad
B_j(\xi)= \max\left\{1-\left|\frac{\xi}{\Delta x}-j\right|,0\right\},
\]
where $B_j$ is the piecewise-linear basis function centered
at $j\Delta x$. This construction interpolates the values
of $g(n^e)-g(1)$ linearly between neighboring grid points.

We evaluate the weights once for each grid using eight point
Gauss-Legendre quadrature on each cell within the sensing
interval. 
Since $h$ is odd, the weights
satisfy $\omega_{-j}=-\omega_j$ and hence
$\sum_{j=-m}^{m}\omega_j=0$. Subtracting $g(1)$ therefore
leaves the weighted sum unchanged in exact arithmetic.
It also ensures that the FFT input is identically zero
at the healed equilibrium. The weighted sum is evaluated by a zero-padded FFT convolution
to prevent periodic wrap around. 

To compute the adhesion flux
derivative, we apply Fourier continuation to the physical
flux values $n_iV_i$, differentiate the resulting Fourier
representation, and retain the values at the physical
interior grid points. The adhesion integral uses the
prescribed collar values, independently of the continuation
used for differentiation.

The spatial discretization, with the prescribed boundary and
exterior collar conditions incorporated, gives a system of differential equations
\[
\frac{d\mathbf{n}}{dt}=\mathbf{F}_h(t,\mathbf{n}),
\]
where $\mathbf{n}(t)$ contains the density values at the
physical interior grid points. The right hand side
$\mathbf{F}_h$ includes the discretized diffusion term,
adhesion flux divergence, and reaction term
$n(1-n)$. We solve this system using the
fourth-order Runge-Kutta method.

\subsection{Spectral thresholds and linear validation}
\label{subsec:num-thresholds}

For the bounded spectral calculation, let
\[
e_j(x)=L^{-1/2}\sin\!\left(\frac{j\pi(x+L)}{2L}\right),
\quad k_j=\frac{j\pi}{2L}.
\]
The Galerkin matrix of the full linearized operator is
\begin{equation}
\label{eq:num-galerkin-matrix}
A^{(m)}_{ij}
=-(\alpha k_j^2+1)\delta_{ij}
+(\lambda-2)\int_{-L}^{L}
e_i'(x)(\mathcal V_{L,\rho}e_j)(x)\,dx,
\quad 1\le i,j\le m.
\end{equation}
This weak formulation incorporates the sensing endpoint
contributions. Individual sine modes are not assumed to be
eigenfunctions of the adhesion operator.

Algorithm~\ref{alg:num-spectral} computes approximations
of the spectral bound $s_L(\rho;\lambda)$ and an associated
eigenfunction using the Galerkin matrix
\eqref{eq:num-galerkin-matrix}.

\begin{algorithm}[ht!]
\caption{Computation of $s_L(\rho;\lambda)$}
\label{alg:num-spectral}
\begin{algorithmic}[1]
\State \textbf{Input:} $L,\rho,\lambda,\alpha,m_0,\varepsilon_s$.
\State \textbf{Output:} Approximate spectral bound $s_L$
and normalized eigenfunction $\phi$.
\State Set $s_{\mathrm{old}}=+\infty$.
\For{$m=m_0,2m_0,4m_0,\dots$}
    \State Assemble $A^{(m)}$ from
    \eqref{eq:num-galerkin-matrix}, using zero exterior
    values and refining quadrature at fixed $m$.
    \State Compute the largest eigenvalue $s$ across both
    symmetry blocks and its unit eigenvector $\mathbf a$.
    \State Check matrix symmetry and the eigenpair residual.
    \If{$|s-s_{\mathrm{old}}|\le\varepsilon_s$}
        \State \Return $s$ and
        $\displaystyle \phi(x)=\sum_{j=1}^{m}a_j e_j(x)$.
    \EndIf
    \State Set $s_{\mathrm{old}}=s$.
\EndFor
\end{algorithmic}
\end{algorithm}

Algorithm~\ref{alg:num-threshold} computes the critical parameters.
For the whole line calculation, the search interval
$0\le k\le k_{\max}$ is chosen with
$k_{\max}>|\lambda-2|H_\rho/\alpha$ at each trial radius.
The estimate
$$\sigma(k)\le-\alpha k^2-1+|\lambda-2|H_\rho k$$
and $\sigma(0)=-1$ ensure that this interval contains a global
maximizer. The wavenumber search is refined to resolve competing
local maxima.

\begin{algorithm}[ht!]
\caption{Critical sensing radii and domain half-length}
\label{alg:num-threshold}
\begin{algorithmic}[1]

\State \textbf{Input:} $\alpha,\lambda$, fixed $L$ or $\rho$
as applicable, search interval $[a,b]$, and tolerance
$\varepsilon_\theta>0$.

\State \textbf{Output:} $\rho_c^\infty(\lambda)$,
$\rho_c^{(L)}(\lambda)$, or $L_c(\rho,\lambda)$,
and the associated critical mode.

\State Choose the threshold and define $G(\theta)$ below.
The diffusion coefficient $\alpha$ remains fixed.
\[
\begin{array}{c|c|c|c}
\text{Threshold} & \text{Fixed} & \text{Varied} & G(\theta)\\ \hline
\rho_c^\infty(\lambda)
& \lambda & \theta=\rho & F(\theta;\lambda)\\[2pt]
\rho_c^{(L)}(\lambda)
& L,\lambda & \theta=\rho & s_L(\theta;\lambda)\\[2pt]
L_c(\rho,\lambda)
& \rho,\lambda & \theta=L & s_\theta(\rho;\lambda)
\end{array}
\]

\State Evaluate $F$ by quadrature and maximization of
$\sigma(k)$ over resolved wavenumbers, including endpoints.
Evaluate $s_L$ using the Galerkin method in
Algorithm~\ref{alg:num-spectral}.

\State Starting from $G(a)<0$, scan $[a,b]$ for the first
transition to $G>0$ when varying $\rho$, or to $G\geq0$
when varying $L$.

\If{the transition is an isolated sign-changing crossing}
    \State Find $\theta_-<\theta_+$ around the first
    detected crossing such that
    \[
        G(\theta_-)<0<G(\theta_+).
    \]
    \State Apply bisection to $G(\theta)=0$ on
    $[\theta_-,\theta_+]$, recomputing $G$ at each iterate,
    until $\theta_+-\theta_-<\varepsilon_\theta$.
    \State Set $\theta_c=(\theta_-+\theta_+)/2$.
\Else
    \State \Return no resolved crossing in $[a,b]$,
    or an unresolved neutral interval or contact.
\EndIf

\State Repeat with finer scanning, quadrature, and
wavenumber or Galerkin resolution until the threshold
changes by less than $\varepsilon_\theta$.
Keep the Galerkin dimension fixed during each bisection solve.

\State \Return the threshold and its critical wavenumber
$k_c$ and wavelength $2\pi/k_c$ for the whole line,
or its normalized eigenfunction $\phi_c$ for the bounded domain.

\end{algorithmic}
\end{algorithm}

We use $\lambda=5,6$ and $L=30,100,200$ for the sensing radius
comparison. The domain-length study starts with
$(\lambda,\rho)=(5,16),(5,18),(6,10),(6,11)$ and refines the
search near $L_c$. Both even and odd perturbations are included.
For the fixed arctangent kernel, the bounded radius scan need
not extend beyond $\rho=2L$, where the operator stops changing.
The spatial structure at the stability threshold is described
by the normalized critical eigenfunction $\phi_c$.

To validate the predicted linear rates, we solve the
model at $L=100$, $\lambda=5$, and
$\rho=16,18$, over $0\le t\le T$, with final time $T=30$
and initial density
\[
n_0(x)=1+\varepsilon
\frac{\phi(x)}{\|\phi\|_{L^\infty}},
\quad \varepsilon=10^{-4},\ 10^{-5},
\]
where $\phi$ is the eigenfunction associated with
$s_L(\rho;\lambda)$. Writing
$a(t)=\langle n(\cdot,t)-1,\phi\rangle$,
we compare $|a(t)|/|a(0)|$ with the predicted exponential $e^{s_Lt}$.
We estimate the observed rate from the slope of a linear
fit to $\log|a(t)|$ against time and compare it with $s_L$,
using times for which the perturbation remains small
and $|a(t)|$ exceeds the numerical error level.

\begin{figure}[t!]
\centering
\includegraphics[width=0.99\linewidth]{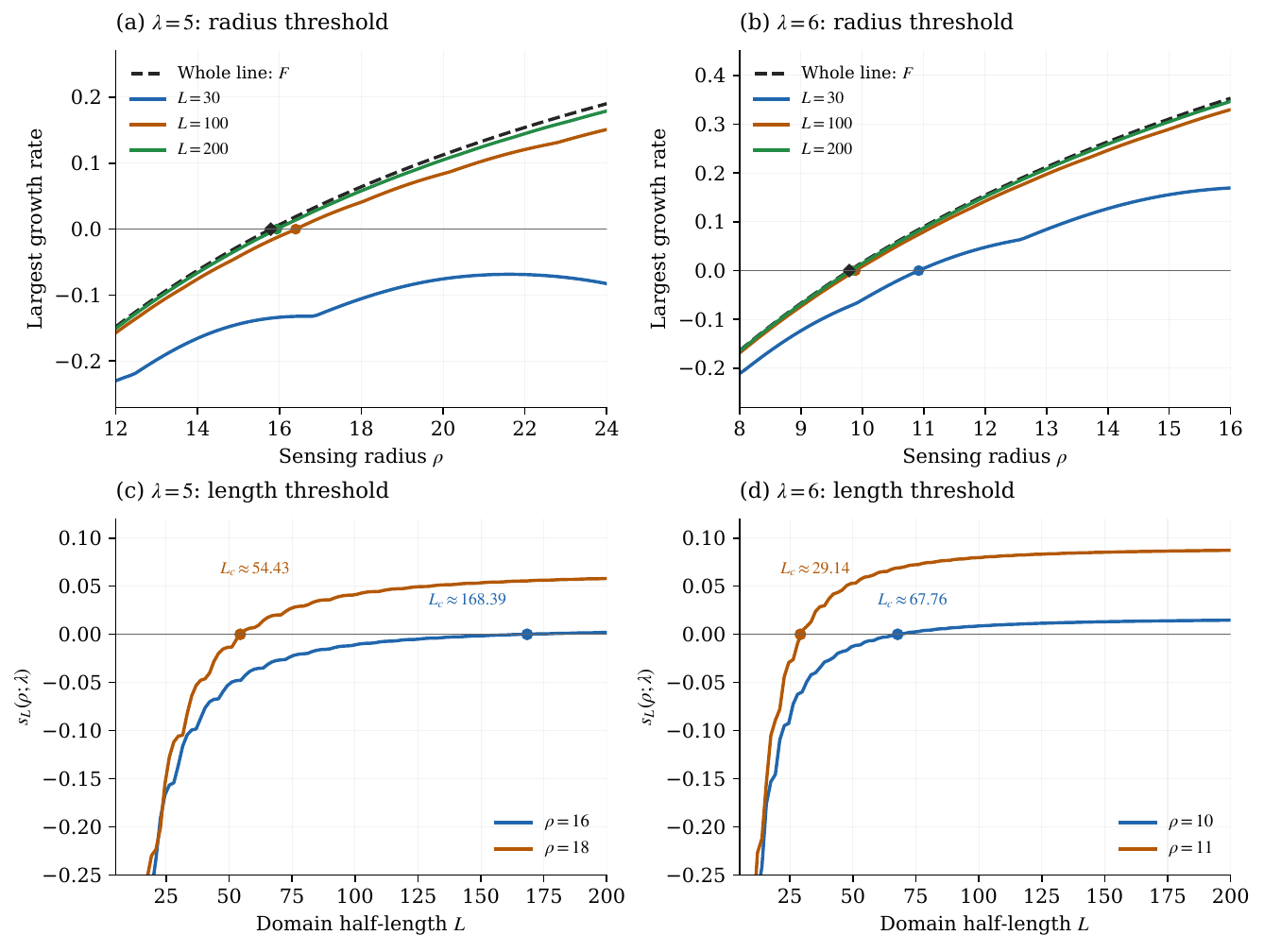}
\caption{Computed spectral thresholds for $\alpha=1$.
Top: the whole line rate $F$ and bounded rates $s_L$ as the
sensing radius varies. Bottom: $s_L$ as the domain half-length
varies. Markers indicate the resolved zero crossings.}
\label{fig:num-linear}
\end{figure}

Figure~\ref{fig:num-linear} gives whole line critical radii
$15.789$ and $9.788$ for $\lambda=5$ and $6$.
The corresponding bounded thresholds at $L=100$ are
$16.397$ and $9.886$.  For $\lambda=5$ and $L=30$, the computed spectral bound satisfies
$s_{30}(\rho;5)\le-0.068$ at the sampled sensing radii
$\rho=0.25,0.50,\dots,60$. Since the perturbation operator does not depend on
$\rho$ once $\rho\ge2L=60$, the same bound holds at every larger radius. At
this domain size the healed state is therefore linearly exponentially stable
for every sensing radius, and $\rho_c^{(30)}(5)=+\infty$ in the sense of
\eqref{eq:bd-critical-radius}. At $\rho=16$, the critical half-length is $L_c\approx168.39$.
Increasing the Galerkin dimension from $384$ to $768$ changes
radius thresholds by less than $4\times10^{-7}$ and length
thresholds by at most $0.0026$.

\begin{figure}[ht!]
\centering
\includegraphics[width=0.95\linewidth]{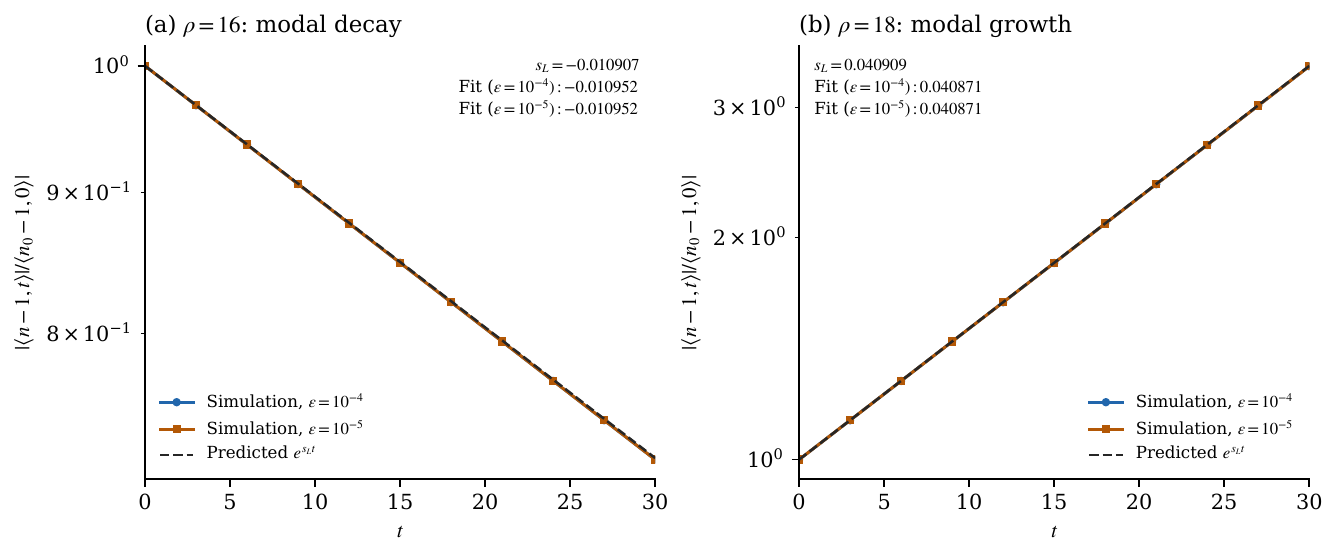}
\caption{Comparison of the predicted linear growth and decay
rates with simulations of the wound healing model at
$L=100$ and $\lambda=5$. Normalized projected amplitudes
are compared with $e^{s_Lt}$ for $\rho=16,18$ and
perturbation sizes $\varepsilon=10^{-4},10^{-5}$.}
\label{fig:num-modal-rates}
\end{figure}

Figure~\ref{fig:num-modal-rates} compares the normalized
amplitudes $|a(t)|/|a(0)|$ with the predicted exponential
$e^{s_Lt}$. The simulations show exponential decay at
$\rho=16$ and growth at $\rho=18$, consistent with the
predicted linear rates for small perturbations.

\subsection{Wound recovery and nonlinear decay}
\label{subsec:num-wounds}

The recovery experiment uses the initial density
\begin{equation}
\label{eq:num-initial-wound}
n_0(x)=1-0.75\exp\!\left[-(0.1x)^6\right],
\quad -L<x<L,
\end{equation}
with $n(\pm L,t)=1$ and $n(x,t)=1$ on the exterior collar
for all $t\ge0$. This smooth, symmetric interior profile has
minimum density $n_0(0)=0.25$ and approaches the equilibrium
density rapidly away from the wound. The profile is held fixed when comparing $L=30,100,200$
at $(\lambda,\rho)=(5,16)$, for which $L_c\approx168.39$.
The endpoint values are set to $1$.

\begin{figure}[t!]
\centering
\includegraphics[width=0.85\linewidth]{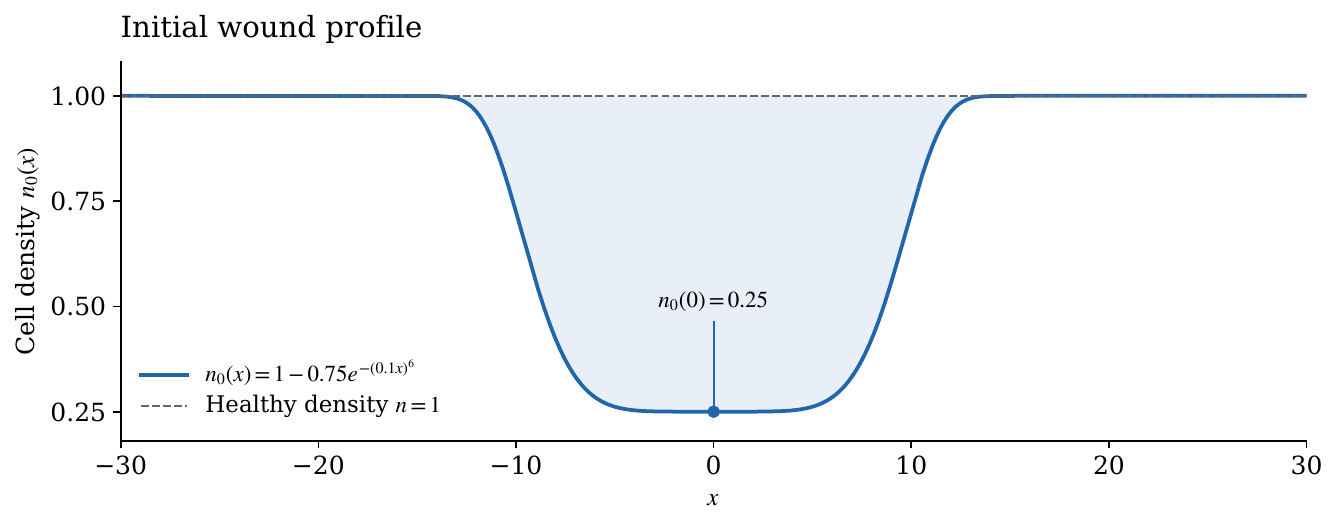}
\caption{Initial recovery profile
$n_0(x)=1-0.75\exp[-(0.1x)^6]$. The initial minimum is $0.25$,
and the prescribed boundary and exterior density is $1$.}
\label{fig:num-initial}
\end{figure}

Figure~\ref{fig:num-initial} shows a central density deficit
of $75\%$, with density close to 1 away from the wound.

We monitor
\begin{equation}
\label{eq:num-diagnostics}
E_L(t)=\frac{\|n(\cdot,t)-1\|_{L^2}}{\sqrt{2L}},
\quad U_L(t)=\|n(\cdot,t)-1\|_{L^\infty},
\end{equation}
and define the observed healing time by
\begin{equation}
\label{eq:num-healing-time}
T_{\mathrm{heal}}^{(T)}
=\inf\left\{t\in[0,T-\tau_h]:
\sup_{s\in[t,T]}U_L(s)\le\varepsilon_h\right\},
\quad \varepsilon_h=10^{-2},\quad\tau_h=20.
\end{equation}
The uniform criterion avoids interpreting spatial averaging
over a larger domain as improved healing. It also requires
recovery to persist for at least $\tau_h$ and through the end
of the observation window.

\begin{figure}[t!]
\centering
\includegraphics[width=\linewidth]{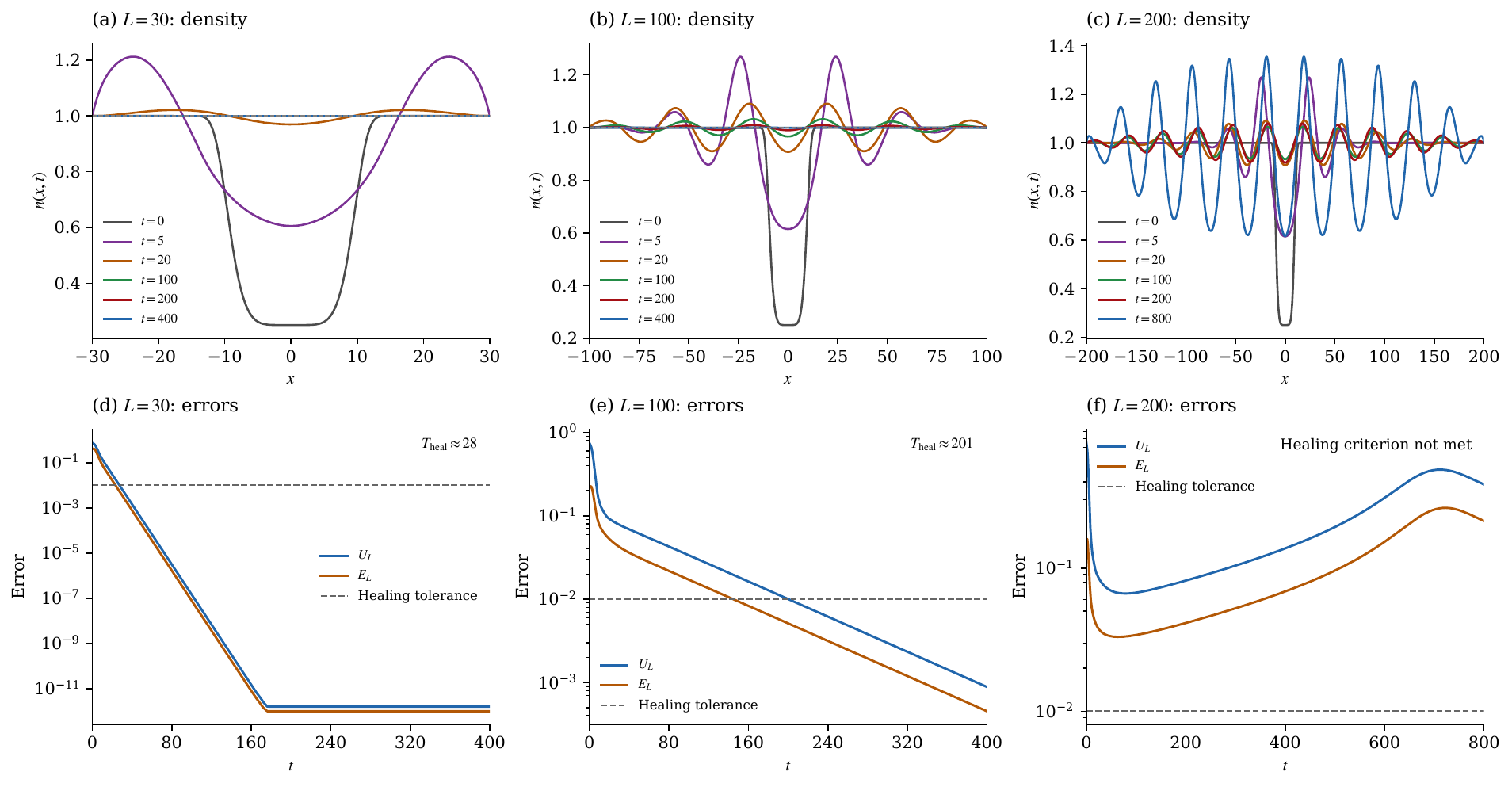}
\caption{Wound healing simulations for
$n_0(x)=1-0.75\exp[-(0.1x)^6]$ and
$(\alpha,\lambda,\rho)=(1,5,16)$.
Top: density profiles at successive times. Bottom:
uniform error $U_L$ and normalized $L^2$ error $E_L$;
the dashed line is the $1\%$ healing tolerance.}
\label{fig:num-wounds}
\end{figure}

Figure~\ref{fig:num-wounds} demonstrates recovery to within
$1\%$ of the healthy density at $t=28$ for $L=30$ and
$t=201$ for $L=100$, sustained through $T=400$.
At $L=100$, the final uniform error is $9.17\times10^{-4}$.
For $L=200$, healing is not observed by $T=800$;
the final uniform error on the finest grid is $0.362$.
These outcomes agree with the signs of $s_L$ in the three
cases, but the finite observation window does not establish
permanent wound arrest.

For the comparison with the sufficient nonlinear estimate, a
separate small wound test uses $n_0(x)=1-d\psi(x/a)$, where
\[
\psi(z)=
\begin{cases}
\exp\!\left(1-\dfrac{1}{1-z^2}\right),&|z|<1,\\
0,&|z|\ge1,
\end{cases}
\]
with
\[
(\alpha,\lambda,\rho)=(1,5,4),\quad
(\delta_0,M)=(0.4,1.2),\quad (a,d)=(0.5,0.5).
\]
Formula~\eqref{eq:bd-Lstab} gives
$L_{\mathrm{stab}}\approx1.79893$. We compare
$L/L_{\mathrm{stab}}=0.6,0.9,1.4$, plotting $E_L(t)/E_L(0)$
against $e^{-\mu_Lt}$ when $\mu_L>0$.
These runs use $\Delta x\approx0.025$,
$\Delta t=10^{-4}$, and output spacing $0.1$, to $T=8$.
The exponential bound is taken from
Theorem~\ref{thm:bd-domain-decay}. The numerically observed density extrema are consistent
with the assumed density bounds, but do not prove that these bounds hold for all time.

\begin{figure}[H]
\centering
\includegraphics[width=0.95\linewidth]{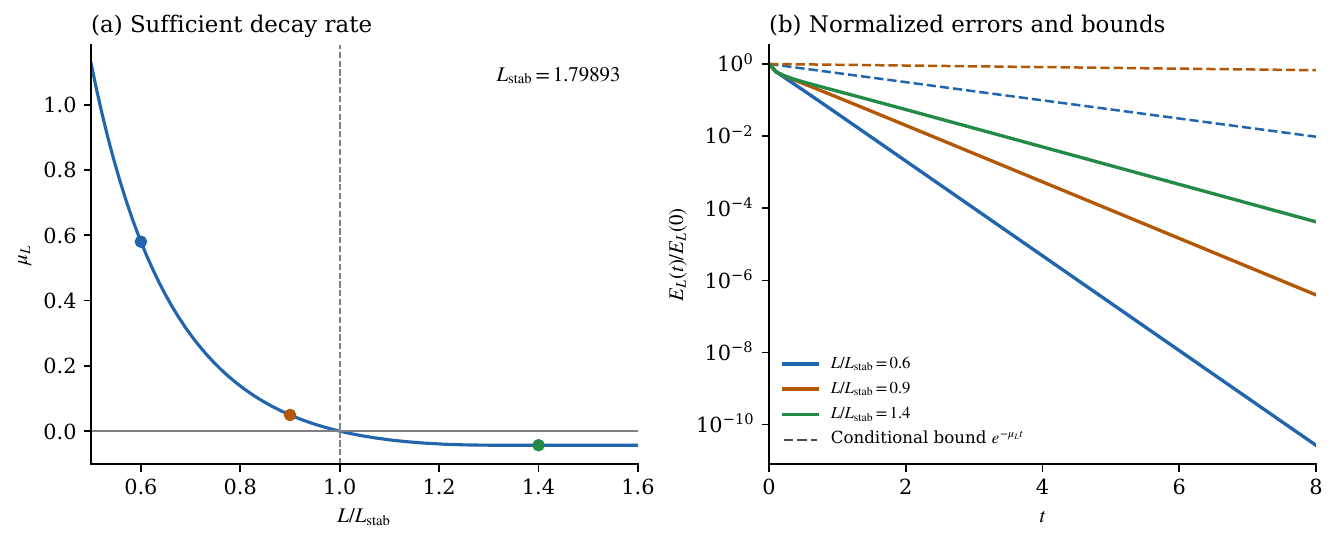}
\caption{Separate small wound test. Left: the sufficient decay
rate $\mu_L$. Right: simulated normalized errors (solid)
and conditional exponential bounds (dashed).
The bounds require global existence and the stated density interval.}
\label{fig:num-decay-bounds}
\end{figure}

In Figure~\ref{fig:num-decay-bounds}, the normalized errors
in both simulations with $L<L_{\mathrm{stab}}$ decay faster
than the corresponding conditional exponential bounds.
The sampled densities remain between $0.5$ and $1$,
consistent with the assumed interval $[0.4,1.2]$.
Decay also occurs at $L=1.4L_{\mathrm{stab}}$, where the
estimate does not guarantee it. Thus the numerical results
illustrate that $L_{\mathrm{stab}}$ is a sufficient bound
rather than a sharp healing threshold.

For the main recovery tests, the last spatial refinement
changes the sampled density histories by at most
$9.22\times10^{-3}$ on $0\le t\le400$;
the independent time-step check at $L=100$ changes them by
less than $5.5\times10^{-10}$.
For $L=200$ at $T=800$, the two finest grids give uniform
errors $0.399$ and $0.362$. Thus the precise late-time profile
remains mesh-sensitive, although neither run satisfies the
healing criterion.
All runs preserve the prescribed boundary values and have
positive sampled densities. These checks quantify numerical
sensitivity; they are not rigorous error bounds.

\section{Conclusion}
\label{sec:conclusion}

We analyze the stability of the healed state $n\equiv1$
in a nonlocal wound healing model incorporating cell diffusion,
proliferation, and adhesion. The whole line dispersion relation
identifies critical sensing radii associated with loss of
linear stability. On bounded intervals, stability also depends on the domain size and on the
prescribed healthy exterior condition that closes the nonlocal operator. The full bounded domain operator determines the
critical sensing radius $\rho_c^{(L)}(\lambda)$ and domain
half-length $L_c(\rho,\lambda)$. The largest bounded growth rate is nondecreasing in $L$ and is bounded above
by the whole line dispersion relation, so whole line linear stability implies
bounded linear stability for this formulation, while the converse may fail. These thresholds cannot,
in general, be obtained by evaluating the whole line
dispersion relation at discrete wavenumbers.

The nonlinear analysis provides a complementary sufficient
condition for healing. Energy estimates and the Poincar\'e
inequality yield a domain dependent decay rate $\mu_L$
and a sufficient half-length $L_{\mathrm{stab}}$.
For $L<L_{\mathrm{stab}}$, solutions converge exponentially
to the healed equilibrium, provided that they exist globally,
remain sufficiently regular, and satisfy the stated density
bounds. Lemma 3.1 gives a sufficient condition on the sensing radius under which the density bounds hold automatically, so that parameters satisfying it only assume global existence and regularity. 

The numerical experiments compare the predicted linear rates with simulations
of the model and quantify the effect of domain size on wound recovery.
Reducing the half-length stabilizes the healed state at adhesion strengths and
sensing radii for which the whole line problem is linearly unstable. At
$\lambda=5$ and $L=30$ no zero crossing of $s_{30}(\cdot;5)$ occurs at the
sampled radii, and since the perturbation operator no longer depends on $\rho$
once $\rho\ge2L$, the healed state at this domain size is linearly
exponentially stable at every sensing radius. They also demonstrate healing
beyond $L_{\mathrm{stab}}$, confirming that the nonlinear
bound is sufficient rather than sharp. The distinction
between $L_c$ and $L_{\mathrm{stab}}$ is therefore essential:
the former describes loss of local linear stability,
whereas the latter provides a conditional guarantee
of exponential healing.

Several questions remain. The simulations show healing
beyond the sufficient bound $L_{\mathrm{stab}}$, motivating
sharper nonlinear criteria. The present analysis fixes the
endpoint and exterior densities at the healed value.
Considering alternative boundary formulations, including
no-flux conditions with adhesive or repulsive boundary
interactions as in~\cite{hillen2020nonlocal}, would allow
us to investigate how these choices affect the equilibria
and their stability thresholds. Extending the analysis to
two spatial dimensions would also allow us to examine
how domain size and shape influence stability, and how
the initial wound geometry affects healing dynamics.

\section*{Code availability}
The code used to produce the numerical results in this paper is publicly
available at \url{https://doi.org/10.5281/zenodo.24356247}.


\bibliographystyle{acm}
\par\vspace{0.3cm}
\bibliography{ref}

\end{document}